\documentclass{article}
\usepackage[utf8]{inputenc}
\usepackage{pgfplots}
\usepackage{amsmath}
\usepackage{amssymb}
\usepackage{mathtools}
\usepackage[]{amsthm}
\usepackage{enumitem}
\usepackage{float}
\usepackage{xcolor}
\usepackage{collectbox}
\usepackage{tikz}
\usepackage{hyperref}
\usetikzlibrary{shapes}
\usetikzlibrary{decorations.markings}
\usetikzlibrary{arrows}

\pgfplotsset{compat=1.16}
\usepackage{geometry}
\usepackage{hhline}

\usepackage{graphicx} 
\begin{document}
\setcounter{page}{1} 
\begin{center}
{\Large\bf Higher diameters of Cayley graphs } 
\end{center}
\medbreak
\centerline{by}
\medbreak
\centerline{Gregory C Magda, Jonathan Rubin, Sabrina Streipert, }

\centerline{Abhiram Kumar,  Cameron Watt}

\centerline{Department of Mathematics}

\centerline{University of Pittsburgh}

\centerline{}

\medbreak
\centerline{and}
\medbreak
\centerline{Gregory P Constantine}

\centerline{School of Computer Science}

\centerline{Georgia Institute of Technology}

\vskip1cm 
\centerline{\bf ABSTRACT}
 
\medbreak\noindent Measures of spread of information in a 
digraph are introduced. Intuitively they assess the 
probability, speed, or number of steps it takes to spread 
information to the entire digraph, thus achieving digraph 
synchrony, by starting from a 
random subset of active vertices. These measures may 
be viewed as 
generalizations of the diameter of a digraph. The paper studies 
these concepts within the context of Cayley digraphs 
associated to finite groups. Explicit formulae for these 
measures are difficult to obtain in general abstract 
groups, and our initial effort is to understand them in 
the case of cyclic groups of prime order. 
With appropriate assumptions on the growth of the 
generating sets, asymptotics show that all the higher 
diameters of random Cayley digraphs are almost surely 
at most 2, as the digraph order goes to infinity.
   
\medbreak\noindent {\em Proposed }
{\em running head:\/} Higher diameters of Cayley digraphs 
\vskip1cm 
\footnoterule\noindent Funded under NIH grant 
RO1-HL-076157 and NSF award ID 2424684 

\noindent gmc@pitt.edu 
  
\newpage
\section*{\Large\bf1.  Motivation and preliminaries} 
\medbreak\noindent The motivating biological problem is to 
turn on all neurons in a brain, or part of a brain, by 
starting with a small subset of active neurons.  We 
view this activity as having {\em local component\/}s, which we 
want to turn on as fast as possible, and {\em global links }
between the local components which serve the purpose 
of efficiently integrating the local components in such a 
way that the entire brain becomes active as quickly as 
possible (or at a controlled prescribed pace).  Further 
detailed information is found in $[14]$ and [13].  Intuitively 
we thus seek to determine those neuronal configurations, 
viewed as abstract networks, that spread the 
information most efficiently to the whole region of 
interest in the brain.  We 
focuss initially on modeling the local components and 
start by making some simplifying assumptions.  The 
basic working hypothesis is that a neuron is activated 
by receiving input from at least $t$ already active 
neurons connected to it.  Initially we make the 
assumption that the underlying digraph that 
connects the neurons is regular; see [1] and
[2].  Since we often want a quick spread to activate 
the whole local area, it is intuitive that the best way of 
doing this is to avoid having short cycles, like triangles 
or squares, since they impede propagation.  If we have $n$ 
neurons, each of (out-)degree $d,$ the emerging 
optimization strategy is that we want to first restrict 
to having a minimal number of triangles, then among 
this subset of regular graphs to seek those that have a 
minimum number of closed walks of length 4 (like 
4-cycles), and proceed sequencially to closed walks of 
higher order.  \medbreak\noindent Imagine for a moment 
that the vertices of the graph (or digraph) are neurons 
and an existing edge transmits information form one 
neuron to another.  We start with a set $S$ of neurons, 
which we call $active,$ and an $activation$ (or $acquiring)$ 
threshold $t,$ which is a natural number.  The spreading 
of neuronal activity is described next.  This is subject 
to some restrictions formulated in terms of {\em Steps}, 
which we now describe.  
\vskip.3cm\noindent
$Step$ 0:  Start with a set $S=S_0$ of vertices of the 
digraph $G$ and a natural number $t.$ We call elements of 
$S$ active vertices.  Imagine that you hold the active 
vertices in your left hand, and the other vertices in 
your right hand.  Color any edge emanating from $S$ red.  
\vskip.3cm\noindent
$Step$ 1:  Acquire vertex $v$, held in your right hand, if $v$ 
has $t$ or more red arrows pointing to it.  [Acquiring a 
vertex is synonomous with the vertex turning active.]  
Move all acquired vertices to your left hand.  Call the 
set of vertices you now hold in your left hand $S_1.$ Color 
all edges emanating from $S_1$ red.  
\vskip.3cm\noindent
The general step is as follows.  We are in posession of 
$S_{i-1}$ with all edges emanating from it colored red.  
\vskip.3cm\noindent
$Step$ $i:$ Acquire vertex $v,$ held in your right hand, if it 
has $t$ or more red arrows pointing to it.  Move all 
acquired vertices to your left hand.  Call the set of 
vertices you now hold in your left hand $S_i.$ Color all 
edges emanating from $S_i$ red.  
\vskip.3cm\noindent
Evidently $S=S_0\subseteq S_1\subseteq\cdots\subseteq S_i\subseteq
\cdots .$ As we keep 
increasing $i$, the following will (obviously) always occur:  
the number of vertices in your right hand becomes 
stationary; that is, $\exists m$ such that, at Step $i,$ for all $
i\geq m$ 
the number of vertices in your right hand remains 
constant.  If your right hand becomes empty for a 
sufficiently large $i$ we say that the digraph is in 
{\em synchrony}.  [You are now holding the whole digraph in 
your left hand -- hence all vertices of the digraph 
became active.]  If the digraph is in a stationary state 
and your right hand is not empty we say that the 
digraph cannot be brought to synchrony when starting 
from subset $S.$ 
\vskip.3cm\noindent
We say that vertex $x$ of digraph $G$ {\em is acquired from }
{\em subset} $S$ if, starting from the active vertices of $S$, 
after a finite number of steps we acquire $x.$ For a given 
starting subset $S$ with $s$ vertices we denote by $d(S,t)$ 
the smallest number of steps that brings the digraph to 
synchrony.  When a digraph cannot be brought to 
synchrony (starting with an incipient set $S=S_0$ and $t)$, 
we convene to write $d(S,t)=\infty$.  \medbreak\noindent Fix 
$t.$ In a digraph $G,$ let $S=S_0$ be a subset with $s$ 
vertices, which we call a $s-$subset.  We introduce the 
following measures of synchrony for $G$.  
\vskip.3cm\noindent
The ratio $p_{st}(G)=$(number of $s-$subsets $S$ that bring $G$ 
to synchrony)/(number of all $s-$subsets) signifies the 
{\em probability\/} (of bringing digraph $G$) {\em to synchrony\/} from a 
randomly chosen $s-$subset.  For given $s$ and $t$, generally 
we are interested in identifying digraphs with large $p_{st}$.  
It might also be observed that there are many instances 
when a digraph has a large $p_{st}$ but the number of steps 
required to attain synchrony are generally quite large, 
which is not so good.  \medbreak\noindent We could tune 
this up by defining another measure $v_{st}(G)$, which we 
call $velocity$ $to$ $synchrony$, as follows:  
\vskip.3cm 
${n\choose s}v_{st}(G)=\sum_Sd(S,t)^{-1}.$ 
\vskip.3cm\noindent
Observe that when $S$ does not induce synchrony, 
$d(S,t)=\infty ,$ and we simply add a zero to the sum.  
Intuitively, velocity $v_{st}$ yields the average speed to the 
synchrony of $G$ across all $s-$subsets.  High values of $v_{st}$ 
are typically good, since the synchrony is then speedily 
restored.  We did not see the concept of velocity to 
synchrony used in the digraph (or network) optimization 
literature so far.  \medbreak\noindent A maxmin analog is 
also very natural to consider.  Define $d_{st}=max_Sd(S,t).$ 
As is easy to see from definitions, $d_{11}$ is just the 
diameter of the digraph.  We call $d_{st}$ the {\em higher }
{\em diameters\/} of digraph $G$ and they are the main subject of 
study in this paper.  

\medbreak\section*{\Large\bf2.  Properties of the measures 
of spread} \medbreak\noindent The following result follows 
immediately from the definitions.  
\medbreak\noindent {\bf Proposition 1} {\em Matrices $(p_{st}),$ $
(v_{st})$ and 
$(d_{st})$ are lower triangular.  In each of these matrices
the entries in each row are decreasing and the entries in each column are increasing.}  \medbreak\noindent The 
word increasing does not mean the increase is strict; 
likewise for decreasing.  Typically we write $p_{st}(H)$ for 
the measure $p_{st}$ of digraph $H,$ and likewise for the 
other measures.  If $G$ is a digraph and $e$ is an edge 
whose endpoints are among the vertices of $G,$ we write 
$H=G\cup e$ for the digraph obtained from $G$ by tossing in 
the extra edge $e.$ (Edge $e$ could be, but typically is not, 
an already existing edge of $G$).  It is again easy to see 
from the definitions (by considering the measure of 
spread in the original digraph without edge $e$ and then 
adding the spread through edge $e$) that the measures we 
introduced are also monotone in any new information 
(that is, edge $e$) added to a digraph; we state this below.  
\medbreak\noindent {\bf Proposition 2} {\em The measures of spread 
have the following monotone properties: $p_{st}(G\cup e)\geq p_{st}(G),$ $d_{st}(G\cup e)\leq d_{st}(G),$ and
$v_{st}(G\cup e)\geq v_{st}(G).$} \medbreak\noindent Recall that a digraph 
is called $regular$ of degree $h$ if each vertex has $h$ edges 
pointing to it, and $h$ edges emanating from it. Let us start with a 
(possibly regular) digraph $G$ with $n$ vertices of degree $h$, 
a threshold $t$ and a subset $S$ of $|S|=s$ active vertices.  
We want to examine what happens when stationarity 
sets in.  By the definition we are in the presence of a 
vertex partition $A$ and $I:=G\backslash A$, where $A(\supseteq S
)$ is the set 
of active vertices acquired from $S$; we assume also that 
$I\neq\emptyset .$ Stationarity means that no vertex in $I$ can be 
acquired in one step from set $A$, and a fortiori cannot be 
acquired from the incipient set $S.$ [For any vertex $x\in G$ 
and any subset $U$ of $G$, let $id_Ux$ denote the in-degree of 
$x$ consisting of all edges that originate in $U$ and point to 
$x;$ analogously $od_Ux$ denotes the out-degree of $x$ 
consisting of all edges that originate at vertex $x$ and 
point to vertices in $U.$] With this notation stationarity 
is equivalent to $\exists A,I$ partition of vertices of $G$ such 
that $id_Ax<t$, $\forall x\in I;$ or, equivalently still (and using the 
regularity of $G)$, $\exists A,I$ partition of vertices of $G$ such 
that $id_Ix>h-t$, $\forall x\in I.$ We now summarize these 
observations in the following result.  
\medbreak\noindent {\bf Theorem 1} 
\vskip.3cm\noindent
{\bf (a)} {\em Let $t$  be the activation threshold, and $s
\geq t.$ Digraph
$G$ cannot be synchronized starting from subset $S$, with
$|S|=s$, if and only if $\exists$ subset $A$ of vertices acquired from $S$, $S\subseteq A\neq G,$ such that $id_A
x<t$, $\forall x\in G\backslash A$.}
\vskip.3cm\noindent
{\bf (b)} {\em Let $t$ be the activation threshold, $h\geq 
t$ and $s\geq t.$ Regular digraph $G$ of degree $h$ cannot be synchronized starting from subset $S$, with $|S|=s$, if and only if $\exists$ subset $A$ of vertices acquired from $S$,$S\subseteq 
A\neq G,$ such that $id_{G\backslash A}x>h-t$, $\forall x\in G\backslash 
A$.}
\vskip.5cm\noindent
\medbreak\noindent\section*{\Large\bf3.  Diameters of 
Cayley graphs} \medbreak\noindent Let $\Gamma$ be a finite group 
and $H$ a subset of $\Gamma .$ Produce a Cayley digraph $Cay(\Gamma 
,H)$ 
of $\Gamma$ and $H$ by defining vertices of $G$ to be elements of 
$\Gamma$ and placing an edge $(x,y)$ from vertex $x$ to vertex $y$ if 
$y=hx$ for some $h\in H.$ [The operations that take place 
are in the group $\Gamma .$] In keeping with tradition, we refer 
to $H$ as a set of generators of $\Gamma$, but {\em no conditions\/} are 
imposed a priori on the subset $H.$ Often the set $H$ is 
understood from the context, in which case we 
abbreviate by simply writing $G$ for $Cay(\Gamma ,H)$.  It is easy 
to see that in a Cayley digraph each vertex has $|H|$ edges 
pointing to it and $|H|$ edges emanating from it; such a 
digraph is called {\em regular of degree} $|H|$ -- we say that 
the in-degree of each vertex is $|H|$ and the out-degree of 
each vertex is $|H|$ as well.  
\vskip.3cm\noindent\ 
Diameters $d_{11}$ of Cayley graphs have been extensively 
studied; cf.  [3], [4] and [5].  Closest to to the 
applications we aim toward 
are examples of large Cayley digraphs and their 
diameters, examples of which appear in [5].  Our main 
interest, however, 
is to initiate the study the higher diameters $d_{st}$.  To our knowleadge 
no such results are available in the literature.  In order 
to do this, for a Cayley digraph $Cay(\Gamma ,H)$ we need to 
start form $s-$subsets and compute $d_{st}$ for fixed $t$$.$ 
Specifically, for a finite group $\Gamma$ fix $H\subseteq\Gamma$ and produce 
$G=Cay(\Gamma ,H)$; whereas the action of $\Gamma$ on the vertices of 
$G$ is transitive (as a right-regular representation), the 
same is not true of the action of $\Gamma$ on the $s-$subsets of 
$G$, for $s>1.$ We need to study this action first.  Let $S$ 
be an $s-$subset of $G$.  We define $Sg$ to be the set 
$\{xg:x\in S\}$.  Group $\Gamma$ acts, therefore, on the $s-$subsets of 
$G$ (or of $\Gamma$) by right translation.  
\medbreak\noindent {\bf Lemma 1} {\em If $\Gamma$ acts by right translation on the s-subsets, then the number of resulting orbits is the coefficient of $y^s$ in the polynomial $P_{\Gamma}(1+y,1+y^2,\ldots ,1+y^n)$.  Here $P_{\Gamma}$ is the cycle index of $\Gamma$ in its right-regular permutation representation on the
group elements of $\Gamma$; $|\Gamma |=n.$} \medbreak\noindent {\bf Proof }
This follows from the well-known orbit-counting theory 
of Redfield and DeBruijn.  For details, see [6], or [7, 
page 225].  \medbreak\noindent The Lemma that follows 
informs us that, in the case of Cayley digraphs, it 
suffices to evolve the spreads simply from 
representatives of the orbits of $s-$subsets.  
\vskip.3cm\noindent
{\bf Lemma 2} {\em For any fixed threshold} $t,$ {\em if s-subsets} $
X$ {\em and }
$Y$ {\em are in the same orbit, then} $d(X,t)=d(Y,t)${\em .  }
\medbreak\noindent {\bf Proof} Write $Y=Xg$ for some $g\in G.$ 
We acquire $z$ from $X$ iff $\exists$ $h_1,\ldots ,h_t\in H$ such that 
$z=h_ix_i$ for $x_i\in X$, $1\leq i\leq t$.  Note now that $z$ is 
acquired from $X$ iff $zg$ is acquired from $Xg=Y$ via the 
same $h_i$ through the images $x_ig\in Xg=Y.$ This bijective 
mapping is being repeated and it holds true at every 
step of the acquiring process.  We conclude that 
$d(X,t)=d(Y,t).$ \medbreak\noindent For any finite group we 
can compute some specific higher diameters of the 
associated Cayley graphs.  \medbreak\noindent {\bf Proposition 3} {\em Let $\Gamma$ be a finite group, $H$ a generating set for $
\Gamma$, and 
$t$ the activating threshold.  In the Cayley graph 
$G:=Cay(\Gamma ,H)$, let $S$ be a set of initially active vertices and $G\backslash S$ be the set of initially inactive vertices.  If 
$t=|H|-k$, then $G$ can be synchronized from $S$ if and only if $G\backslash S$ contains no subdigraph $K$ with $
id_Kx\geq k+1,$ 
$\forall x\in K$.}  \medbreak\noindent {\bf Proof} The result follows from 
Theorem 1, but we offer a self-contained proof.  Let 
$I:=G\backslash S.$ Assume that $I$ contains such a subdigraph $K$.  
We show that $G$ cannot synchronize.  Evidently 
$id_Ix\geq k+1$ for all vertices $x\in K.$ Thus at the initial 
acquiring step $id_Sx\leq |H|-(k+1)<|H|-k=t,$ for all $x\in K.$ 
Therefore no vertices in $K$ will activate at the first 
step.  Iteratively, vertices of $K$ will never activate and 
$G$ will fail to synchronize.  
\vskip.3cm\noindent
We now assume that $G$ does not synchronize and force 
the existence of a subdigraph $K$ in $I$ having the stated 
properties.  The assumption implies that we reach 
stationarity with a bipartition into active vertices $A$ and 
inactive vertices $J$ (with $J$ nonempty).  Clearly $S\subseteq A$ 
and $J\subseteq I.$ Define $K$ to be subdigraph $J.$ Let $x\in K.$ Then 
$id_Kx\geq k+1,$ else $x$ would have been acquired from $A,$ 
since $t=|H|-k.$ It is now seen that $K$ is a subdigraph of 
$I=G\backslash S$ with the desired properties.  \medbreak\noindent
By a cycle in a digraph we always mean a directed 
cycle.  The case of Proposition 3 with $k=0$ is 
highlighted next, since it allows computation of some 
higher diameters.  \medbreak\noindent {\bf Corollary 1} {\em Let $
\Gamma$ be a finite group, $H$ a generating set for $\Gamma$, and $
t$ the activating threshold.  In the Cayley graph $G:=Cay(\Gamma ,
H)$, let $S$ be a set of initially active vertices and $G\backslash 
S$ be the set of initially inactive vertices.  If $t=|H|$, then $
G$ can be synchronized from $S$ if and only if $G\backslash 
S$ is cycle-free.}  \medbreak\noindent {\bf Proof} Taking $k=0$ in 
Proposition 3 yields the result, but we shall examine the 
nature of the subdigraph $K$ more closely in this case.  
Specifically, we assume that $G$ does not synchronize and 
force the existence of a cycle in $I.$ The assumption 
implies that we reach stationarity with a bipartition 
into active vertices $A$ and inactive vertices $J,$ and $J$ 
contains an inactive vertex $c_1.$ Clearly $S\subseteq A$ and $J\subseteq 
I.$ 
Since $t=|H|$ and $c_1\in J$ it follows that there exists 
$c_2\in J$ and $c_2$ points to $c_1.$ We now start with $c_2$ and 
iterate the process.  This yields a set of vertice ($c_1,$ $c_2,$ 
$c_3,$ \ldots ) with $c_{i+1}$ pointing to $c_i$ and with all $c_
j\in J.$ 
Since $J$ is finite this forces one of the $c_i$ to be 
revisited, thus generating a cycle $C$ in $J\subseteq I.$ 
\medbreak\noindent We are now in position to compute 
some higher diameters of Cayley graphs.  Let $G$ be the 
Cayley graph $Cay(\Gamma ,H)$, with $|G|=n$ and $|H|=h.$ Denote 
by $m$ the minimum cycle length among all cycles in $G.$ 
Scalar $m$ is called the girth of $G.$ We demonstrate the 
following result.  \medbreak\noindent{\bf Theorem 2}
{\em Let $G=Cay(\Gamma ,H)$ be a Cayley graph having 
girth $m$.  If the activating threshold $t$ is equal to 
$h=|H|$, then the diameter $d_{st}(G)=n-s$ whenever 
$s>n-m,$ and is infinite otherwise.} 
\medbreak\noindent {\bf Proof} For $s>n-m$, the initial inactive 
set will not contain a cycle so by Corollary 1 digraph $G$ 
will synchronize for all starting conditions with $s$ 
vertices active.  Since there are $n-s$ initially inactive 
vertices, the worst case activation among the initially 
active sets of size $s$ is when the digraph activates in 
exactly $n-s$ steps (i.e.,  the worst case activation is by 
activating vertices one-at-a-time).  In particular, for 
$s>n-m$, we have that $d_{st}(G)\leq n-s$ (We remark that 
this inequality holds for $s$ large enough to guarantee 
that the graph will eventually fully activate).  We will 
show that, in fact, $d_{st}(G)=n-s$.  We assert that this 
value is achieved if we take our initial inactive set to 
be $\{g_1,g_2,...,g_{s-(n-m)}\}$ and all other vertices active (that 
is, take all vertices not in the minimal cycle to be 
active, and then take the first $s-(n-m)$ vertices in the 
cycle to be active).  Then the set of initially inactive 
vertices is $U=\{g_{s-(n-m)+1},...,g_m\}$.  Consider the vertex 
$g_{s-(n-m)+1}$.  We will show that this vertex receives an 
input of $h$ at the first time step.  If we assume for the 
sake of contradiction that $g_{s-(n-m)+1}$ receives an input 
of less than $h$, then there must be some vertex 
$g_k\in U\backslash \{g_{s-(n-m)+1}\}$ that has an edge to $g_{s-
(n-m)+1}$ since 
all other vertices are active.  However, this would form 
a cycle with length smaller than $m$ which contradicts $m$ 
being the minimal cycle length.  Thus, $g_{s-(n-m)+1}$ will 
activate at the first time step.  Now consider some 
vertex $g_k\in U\backslash \{g_{s-(n-m)+1}\}$.  We will show that $
g_k$ does 
not activate at the first time step.  Since 
$k>s-(n-m)+1$, vertex $g_{k-1}$ is an initially inactive.  
Thus, $g_k$ can have at most an input of $h-1$ since $g_{k-1}$ is 
inactive and $g_k$ has an in-degree of $h$.  Thus, $g_k$ will not 
activate for $k\in \{s-(n-m)+1,...,m\}$.  Iteratively, we 
conclude that the remainder of the cycle will activate 
one vertex at a time hence, for $s>n-m$, we obtain 
$d_{st}(G)=n-s$.  
\medbreak\noindent
{\bf Examples} Consider $G=Cay(C_7,H)$ with $H=\{a,a^2,a^3\}.$ The 
girth is 3 as is seen in the cycle $(1,a^3,a^6).$ Theorem 2 
yields $d_{s3}(G)=7-s$ for $s>4,$ and $\infty$ otherwise; this 
can be verified directly. 
\vskip.3cm\noindent
If we consider $s=5,$ then $d_{52}(G)=2;$ this diameter value 
is obtained from the initial active set $\{1,a,a^2,a^4,a^5\}.$
\vskip.3cm\noindent
Lastly, we examine Theorem 2 in the case 
$G=Cay(C_{11},H)$ with $H=\{a,a^6\}.$ Here the girth is 6, 
which is seen in the cycle $(1,a^6,a^7,a^8,a^9,a^{10}$). By 
Theorem 2 we have $d_{s2}(G)=11-s$ for $s>5.$ For 
instance, when $s=6$ we obtain $d_{62}(G)=11-6=5$; this 
value is attained from the initial active set 
$\{1,a,a^2,a^3,a^4,a^5\}.$
\medbreak\noindent\section*{\Large\bf3.  The higher 
diameters of groups of prime order} \medbreak\noindent In 
general it is difficult to find explicit formulae for the 
diameter of sufficiently complex classes of groups, such 
as finite groups of Lie type or symmetric groups.  
Asymptotic approximations exist but the proofs are quite 
intricate and typically deal with undirected graphs (cf.  [4], [8], [9]).  Computing higher diameters does 
not, in any way, simplify the process.  We therefore 
tread modestly by first examining the cyclic groups $C_p$ 
for $p$ prime.  The first thing we observe is that we can 
select (Dyck) paths as orbit representatives for the 
action of $C_p$ on the $k-$subsets of vertices, as explained 
below.  \medbreak\noindent Start with the cyclic group 
$\Gamma =C_n$.  An upward-and-right moving path on the integer 
lattice starting at (0,0) and ending at $(n-k,k)$ that 
touches or stays above the line joining (0,0) and ($n-k,k)$ 
is, for simplicity, just called a $path.$ Note first that any 
path may be viewed as a $k-$subset (and, equivalently, as 
a 0,1-binary sequence with $k$ ones), by listing the indices 
of the upward moves as elements of the set.  To be 
specific, an upward move is marked by 1 and a move to 
the right by 0.  We draw attention to the following 
useful observation.  \medbreak\noindent {\bf Lemma 3} {\em The 
$C_n$-orbit of any binary sequence of length $n$ and weight 
$k$ contains at least one path.  It contains exactly one 
path if $n$ and $k$ are coprime.}  \medbreak\noindent {\bf Proof}
This is best seen as follows.  The line joining $(0,0)$ and 
$(n-k,k)$ has slope $\frac k{n-k}$.  For any 0-1 sequence 
$s=(s_i:1\leq i\leq n)$ of length $n,$ replace the 0s by $\frac {
-k}{n-k}$ and 
leave the 1s as they are; call the new sequence $s'=(s_i').$ 
Calculate the partial sums $p_i$ of $s'$ by setting $p_1=s_1',$ 
and $p_i=p_{i-1}+s_i'$, $2\leq i\leq n.$ Let $i^{*}$ be an index for which 
$p_{i^{*}}$ is a minimum.  Apply a cyclic rotation to the 
original sequence $s$ that places in position 1 the index 
$i^{*}+1.$ It is evident from this construction, by the choice 
of $i^{*}$, that the resulting sequence is in the $C_n$-orbit of $
s$ 
and that it corresponds to a path.  When $k$ and $n$ are 
coprime the index $i^{*}$ is unique; else the path would 
touch the diagonal line at $(h,v)$ with $(h,v)\neq (0,0)$ or 
$(n-k,k)$.  Similar right triangles now yield $\frac vh=\frac k{n
-k},$ or 
$vn=k(v+h).$ Since $k$ and $n$ are coprime this forces $k$ to 
divide $v,$ but this is not possible because $v<k.$ This 
ends the proof.  
\vskip.3cm\noindent For 
instance, if $s=(0,1,0,0,1,1,0),$ then its orbit is 
represented by the shifted sequence (1, 1, 0, 0, 1, 0, 0) 
which corresponds to a path.  Note that $s$ itself does 
not correspond to a path.\medbreak\noindent A consequence 
of the Lemma 3 for prime $n$ may be summarized as 
follows.  \medbreak\noindent {\bf Lemma 4} {\em Let $C_p$ act on 
$k$-subsets of itself by right translation.  Then the 
resulting orbits are in bijection with paths either touching or above the line joining $(0,0)$ to $(p-k,k
).$} 
\medbreak\noindent Examining Lemma 1 and noting that the 
cycle index of $C_p$, acting on its elements by right 
translation, is $P_{C_p}=\frac 1p[(p-1)x_p+x_1^p],$ we conclude that 
there are exactly $\frac 1p{p\choose k}$ orbits on $k-$subsets.  [In this 
case it is, of course, also easy to count the number of 
such orbits directly.]  As an example, we use the paths 
as orbit representatives for an exhaustive computation 
of all higher diameters for the group $C_7.$ The results 
are summarized in Table 1.  \medbreak\noindent
\centerline{{\bf Table 1} Smallest $d_{st}$ diameters for 
$\textrm{Cay}(C_7,H)$.}

\begin{table}[H]
    \centering
    \begin{tabular}{|c|c|c c|c c c|c c c c|c c c c c|}
        \hline
         & $|H| = 1$ & \multicolumn{2}{|c|}{$|H|=2$}  & \multicolumn{3}{|c|}{$|H|=3$}
         \\ \hline
         & $t=1$ & $t=1$ & $t=2$ & $t=1$ & $t=2$ & $t=3$
         \\ \hline
         $s=1$ & 6 & 3 & & 2 & & \\ 
         $s=2$ & 5 & 2 & $\infty$ & 2 & 3 & \\
         $s=3$ & 4 & 2 & $\infty$ & 2 & 2 & $\infty$\\
         $s=4$ & 3 & 2 & 3 & 1 & 2 & $\infty$ \\
         $s=5$ & 2 & 1 & 2 & 1 & 1 & 2 \\
         $s=6$ & 1 & 1 & 1 & 1 & 1 & 1 \\ \hline
    \end{tabular}
    \label{|H| = 1,2,3}
\end{table}

\begin{table}[H] 
\centering 
\begin{tabular}{|c|c 
c c c|c c c c c|}  \hline & \multicolumn{4}{|c|}{$|H|=4$}  & 
\multicolumn{5}{|c|}{$|H|=5$}  \\ \hline & $t=1$ & $t=2$ & 
$t=3$ & $t=4$ & $t=1$ & $t=2$ & $t=3$ & $t=4$ & $t=5$  \\ 
\hline $s=1$ & 2 & & & & 2 & & & & \\ $s=2$ & 2 & 2 & & & 1 
& 2 & & & \\ $s=3$ & 1 & 2 & $\infty$ & & 1 & 1 & 2 & & \\ $s=4$ 
& 1 & 1 & 2 & $\infty$ & 1 & 1 & 1 & $\infty$ & \\ $s=5$ & 1 & 1 & 1 & 
$\infty$ & 1 & 1 & 1 & 1 & $\infty$ \\ $s=6$ & 1 & 1 & 1 & 1 & 1 & 1 & 1 
& 1 & 1 \\ \hline 
\end{tabular}
\label{|H| = 4,5} 
\end{table}
\vskip.3cm\noindent
For $|H|=6$ we have the complete digraph on 7 vertices, 
hence for $s\geq t$, all $d_{st}$ diameters are 1.  
\medbreak\noindent We specify the choices of $H$ that yield 
the smallest $d_{st}$ values in Table 1.  In many cases the 
$d_{st}$ entry is invariant to the choice of $H.$ The cases 
where it is not are as follows.  Write $|H|=h.$ Then 
$(h,t,s)=(2,1,2)$ has $H=\{a,a^3\};$ $(h,t,s)=(2,2,4)$ has 
$H=\{a,a^2\}$; $(h,t,s)=(2,2,5)$ has either $H=\{a,a^2\}$ or 
$H=\{a,a^3\}.$ Also, $(h,t,s)=(3,2,2)$ has $H+\{a,a^2,a^4\}$; 
$(h,t,s)=(3,2,3)$ has either $H=\{a,a^2,a^3\}$ or $H=\{a,a^2,a^4\}$; 
$(h,t,s)=(3,3,5)$ has either $H=\{a,a^2,a^3\}$ or $H=\{a,a^2,a^4\}$.  
Lastly, $(h,t,s)=(4,2,2)$ has $H=\{a,a^2,a^3,a^5\}$; and 
$(h,t,s)=(4,3,4)$ has either $H=\{a,a^2,a^3,a^4\}$ or 
$H=\{a,a^2,a^3,a^5\}.$ Lemmas 2 and 4 were used in 
computing all entries in Table 1 and the corresponding 
generating sets we just listed.  For instance, the 
$\frac 17{7\choose 3}=5$ paths representing the orbits for $p=7,$ $
s=3$ 
are 1110000, 1101000, 1100100, 1011000 and 1010100.  In 
calculations involving entries in Table 1, we used these 
paths for selecting the initially active sets $S.$ 
\medbreak\noindent We may apply Theorem 2 to obtain 
higher diameter values for the group $C_p.$ We nned a 
formula for girth of the Cayley graphs that arise in this 
case. Initially we examine the case $G=Cay(C_p,H)$ where 
$H=\{a,a^k\}.$ The minimum cycle length is easily seen to 
be written as $m=min_{z,w}\{z+w:z+wk=0$ mod $p$$\}$, where 
$z$and $w$ 
are nonnegative integers. We have that $z+wk=np$ for 
some $n\in N.$ Then $z=np-wk,$ and 
$z+w=(np-wk)+w=np-w(k-1).$ In this rewriting,
\vskip.3cm
\centerline{ }
\[m=min_{n,w}\{np-w(k-1):np-wk>0\},\]
\vskip.3cm\noindent
with $n\in N$ and $w\in N\cup 0.$ For a fixed $w,$ the value of $
n$ that attains this minimum 
is given by the minimumvalue of $n$ satisfying 
$np-wk>0.$ For fixed $w$, we may thus take $n=\lceil\frac {wk}p\rceil$, 
where $\lceil x\rceil$ denotes the least integer greater than or 
equal to $x.$ We thus obtain
\vskip.3cm
\[m=min_{w\in \{1,\ldots ,p\}}[p\lceil\frac {wk}p\rceil -w(k-1)].\]
\vskip.3cm\noindent
The argument presented above generalizes and allows us 
to establish the following result.
\medbreak\noindent
{\bf Proposition 4} {\em The girth of $G=Cay(C_p,H)$ with 
$H=\{a,a^{k_1},\ldots ,a^{k_q}\}$ is
\vskip.3cm
\[m=min_{w_1,\ldots ,w_q\in \{1,\ldots ,p\}}[p\lceil\frac 1p\sum_{
i=1}^qw_ik_i\rceil -\sum_{i=1}^qw_i(k_i-1)].\]
}
\vskip.3cm\noindent
Higher diameter values of Cayley digraphs of $C_p$ are 
now explicitely following from Propositions 3 and 
Corollary 1 and Theorem 2.
\medbreak\noindent\section*{\Large\bf4.  Asymptotics on 
the higer diameters of Cayley digraphs} 
\medbreak\noindent Let $\Gamma$ be a group of order $n$ and $H\subseteq
\Gamma$ 
a generating set; $|H|=h.$ Write $p:=\frac hn.$ For a fixed scalar $
c$, $0<c<1,$ 
we assume that $n$ and $h$ can 
vary, but force $p=\frac hn\geq c.$ Let $G$ be the 
Cayley digraph obtained from $\Gamma$ and $H.$ All subsets $H$ with 
$|H|=h$ are assumed equally likely.  Starting from any 
vertex $x$ of the Cayley digraph $G$ the probability that for 
any vertex $y$ the edge $(x,y)$ exists is ${{n-1}\choose {h-1}}{n\choose 
h}^{-1}=\frac hn=p,$ 
for all $n.$ [It is not a serious issue, but we allowed loops in the digraph $
G.$] Since $p$ is free of choices of $n$ and $h,$ 
we prove the following asymptotic result. Although the 
context of randomness differs, the main argument is 
similar to that used for ordinary diameters in [10].  
\medbreak\noindent {\bf Theorem 3} {\em {\rm A Cayley digraph with $
n$ }
vertices and in-degree $h$, having the property that $
p=\frac hn\geq c$
$(0<c<1)$ irrespective of $n$ and h, has almost surely 
diameter $d_{s1}\leq 2$ as $n\rightarrow\infty,$ for all $
s\geq 1.$}
\medbreak\noindent {\bf Proof} Fix two distinct vertices $x$ and 
y.  The probability that we have edges $(x,z)$ and $(z,y),$ $z$ 
different from $x$ and $y$, is $p^2.$ That no such directed path 
through $z$ exists carries chance $1-p^2.$ That there is no 
such path at all between $x$ and $y$ is $(1-p^2)^{n-2},$ since 
there are $n-2$ independent choices for $z.$ As there are 
${n\choose 2}$ choices for the starting vertices $x$ and $y,$ we 
conclude that the probability of no two distinct vertices 
being at distance exactly 2 is ${n\choose 2}(1-p^2)^{n-2}.$ Since $
p\geq c,$ 
$0<c<1,$ we now see that ${n\choose 2}(1-p^2)^{n-2}\geq{n\choose 
2}(1-c^2)^{n-2}$
goes to 0 as $n\rightarrow\infty .$ We conclude that the diameter of the 
Cayley graph is almost surely 2, as $n\rightarrow\infty .$ Use of 
Proposition 1 now ends the proof.  \medbreak\noindent The 
result below is presented in the context of random 
graphs but, just like Theorem 3, it can be also 
formulated for ''random'' Cayley graphs. In this sense it 
generalizes the results in [11] and [12] to the higher 
diameters of random digraphs.  
\medbreak\noindent {\bf Theorem 4} {\em A random digraph in
which each edge is randomly and independently
generated with probability p, $0<p<1$ $($ regardless of
the order of the digraph$)$, has almost surely diameter
$d_{st}\leq 2$ as $n\rightarrow\infty,$ for all fixed $
s\geq t\geq 1.$} 
\medbreak\noindent {\bf Proof} Let $G$ be a random graph with $n$ 
vertices.  Fix a subset $S$ of $s$ vertices, and fix also a 
vertex $v\notin S.$ We compute the probability of acquiring $v$ 
from $S$ in exactly two steps.  To acquire $v$ from $S$ in 
two steps, we must first acquire a vertex $x$, $x\notin S$ and 
$x\neq v,$ and then acquire $v$ from $S\cup \{x\}.$ Acquiring $x$ from 
$S$ carries probability $r_1=\sum_{k\geq t}{s\choose k}p^k\geq p^
t.$ Likewise, the 
probability of acquiring $v$ from $S\cup \{x\}$ is 
$r_2=\sum_{j\geq t-1}p{s\choose j}p^j\geq p^t$; note that in the sum we 
necessarily pick the edge from $x$ to $v$ with probability 
$p.$ The chance of a{\em c\/}quiring $v$ in two steps, as indicated, 
is $r_1r_2\geq p^{2t}.$ To not acquire $v$ in such manner, via $x$, 
has chance $1-r_1r_2\leq 1-p^{2t}.$ There are $n-s-1$ choices for 
the intermediate vertex $x$, all producing independent 
events.  Hence not reaching $v$ form $S$ in two steps 
carries the chance $(1-r_1r_2)^{n-s-1}\leq (1-p^{2t})^{n-s-1}.$ It follows 
that the probability of not reaching from any $s-$subset 
to any vertex not in that subset in two steps is 
${n\choose {s+1}}(1-r_1r_2)^{n-s-1}\leq{n\choose {s+1}}(1-p^{2t})^{
n-s-1}.$ Since $p$ does not 
depend on $n$, and $s$ and $t$ are fixed, it is clear that we 
compare exponential and polynomial growths.  The 
conclusion is that the probability of not reaching from 
any $s-$subset to any vertex not in that subset in two 
steps goes to zero as $n\rightarrow\infty .$ The truth of the enunciated statement now follows.  
\vskip1cm 
\centerline{{\bf Acknowledgement}}
\vskip.3cm\noindent\ 
We are grateful to the National Science Foundation for 
sponsoring this work under the NSF Emerging 
Mathematics in Biology grant 2424684.  
\vskip1cm 
\newpage

\centerline{{\bf REFERENCES}}
\begin{enumerate}
\item Chartrand G, Jordon H, Vatter V, Zhang P {\em Graphs }
{\em and digraphs}, Chapman and Hall, 2024

\item L\"oh, C {\em Geometric group theory}, Springer, 2017

\item Babai L, Seress A (1988) On the diameters of Cayley 
graphs of the symmetric group, {\em J Comb Theory, Series }
{\em A}, 49, 175-179

\item Helfgott H, Seress A (2014) On the diameter of 
permutation groups, {\em Annals of Mathematics,\/} 179, 611-658

\item Erskine G, Tuite J (2018) Large Cayley graphs of 
small diameter, {\em Disc Applied Math}, 250, 202-214

\item Redfield, J H (1927) The theory of group-reduced 
distributions, {\em Amer. J. Math.,\/} {\bf 49}, 433-455

\item Constantine, G {\em Combinatorial theory and statistical }
{\em design}, Wiley, New York, 1987 

\item Tao, T Expansion in finite simple groups of Lie 
type, Graduate studies in mathematics, vol 164, {\em Amer }
{\em Math Society}, 2015

\item Liebeck M, Shalev A (2001) Diameters of finite 
simple groups:sharp bounds and applications, {\em Annals of }
{\em Mathematics}, 154, 383-406

\item Philips T K, Towsley D F, Wolf J K (1990) On the 
diameter of a class of random graphs, {\em IEEE }
{\em Transactions on Information Theory,\/} 36, 285-288

\item Meng J, Huang Q (1994) Almost all Cayley graphs 
have diameter 2, {\em Discrete Mathematics,\/} 178, 285-288

\item Erskine G (2015) Diameter 2 Cayley graphs of 
dihedral groups, {\em Discrete Mathematics,\/} 338, 1022-1024

\item Bohnen N,  Prabesh K,  Koeppe R,  Catasus C,  
Frey K,  Scott P,  Constantine G,  Albin R, M\"uller 
M (2021) Regional cerebral cholinergic nerve terminal 
integrity and cardinal motor features in Parkinson's 
disease, {\em Brain communications,\/} vol 3, issue 2, fcab109 

\item Bear M,  Connors B,  Paradiso M  {\em Neuroscience:  }
{\em Exploring the brain,\/} Fourth edition, Jones and Bartlett 
Learning, Burlington, MA, 2016 

\end{enumerate}

\end{document}